  \documentclass[11pt]{article}

\usepackage{amsmath, amssymb, amsthm, latexsym, bbm}

 \usepackage{color}

 \allowdisplaybreaks
\usepackage{hyperref}

 \usepackage{graphicx}

\newtheorem{thm}{Theorem}[section]
\newtheorem{lemma}[thm]{Lemma}
\newtheorem{prop}[thm]{Proposition}
\newtheorem{defn}[thm]{Definition}

\newcommand{\be}{\begin{equation}}
\newcommand{\ee}{\end{equation}}
\newcommand{\bea}{\begin{eqnarray}}
\newcommand{\eea}{\end{eqnarray}}
\newcommand{\beaa}{\begin{eqnarray*}}
\newcommand{\eeaa}{\end{eqnarray*}}
\newcommand{\bei}{\begin{itemize}}
\newcommand{\eei}{\end{itemize}}
\newcommand{\bee}{\begin{enumerate}}
\newcommand{\eee}{\end{enumerate}}
\newcommand{\bi}{\begin{itemize}}
\newcommand{\ei}{\end{itemize}}
\newcommand{\beq}{\begin{eqnarray*}}
\newcommand{\eeq}{\end{eqnarray*}}
\newcommand{\beqn}{\begin{eqnarray}}
\newcommand{\eeqn}{\end{eqnarray}}

\def\1{{\mathbbm 1}}

\newcommand{\pf}{\noindent {\bf Proof.\ }}

\newcommand{\ignore}[1]{}{}

\def\<{\langle}
\def\>{\rangle}
\def\loc{{\rm loc}}

\def\qed{{\hfill $\Box$ \bigskip}}

\def\P{{\mathbb P}}

\def\E{{\mathbb E}\,}

\def\R{{\mathbb R}}

\def\K{\mathbb K}

\newcommand{\LL}{{\mathcal L}}

\newcommand{\sL}{{\mathcal L}}

\def\eps{\varepsilon}

\numberwithin{equation}{section}

\begin{document}

 \title{\bf Dirichlet  problem for diffusions with  jumps}

\author{{\bf Zhen-Qing Chen} \quad and \quad {\bf Jun Peng}}

\date{\today}
\maketitle

\begin{abstract}

In this paper, we study   Dirichlet problem for non-local operator on bounded domains in $\R^d$
 $$
\mathcal{L}u =  {\rm div}(A(x) \nabla (x)) +  b(x) \cdot \nabla u(x)
  + \int_{\R^d}    (u(y)-u(x)  ) J(x, dy) ,
$$
where $A(x)=(a_{ij}(x))_{1\leq i,j\leq d}$ is a measurable $d\times d$ matrix-valued function on $\R^d$ that is uniformly elliptic and bounded,
$b$ is an $\R^d$-valued function so that $|b|^2$ is in some Kato class $\K_d$,
for each $x\in \R^d$, $J(x, dy)$ is a finite measure on $\R^d$ so that $x\mapsto J(x, \R^d)$  is in  the Kato class $\K_d$.
We show  there is a unique Feller process $X$ having strong Feller property associated with $\sL$,
which can be obtained from the diffusion process  having  generator
 $ {\rm div}(A(x) \nabla (x)) +  b(x) \cdot \nabla u(x) $ through  redistribution.
We  further show that for any  bounded connected open subset $D\subset \R^d$  that is  regular with respect to
the Laplace operator $\Delta$ and for any bounded continuous function $\varphi $ on $D^c$,
the Dirichlet problem $\sL u=0$ in $D$ with $u=\varphi$ on $D^c$ has a unique bounded continuous weak solution on $\R^d$.
This unique weak solution can be represented in terms of the  Feller process
associated with $\sL$.

\bigskip

 
\medskip


\end{abstract}

\section{Introduction}

There is a rich interplay between probability and analysis.
In a pioneering paper  \cite{Ka},  Kakutani showed that for any Jordon domain $D$ in $\R^2$ and any bounded continuous function
$\varphi$ on $\partial D$, $u(x):=\E_x [ \varphi (X_{\tau_D})]$ is a harmonic function in $D$ that is continuous up to the boundary $\partial D$
with   boundary value $\varphi$; that is, $u(x)=\E_x [ \varphi (X_{\tau_D})] $ solves the Dirichlet boundary value problem
 $$
   \Delta u =0  \  \hbox{  in } D  \qquad \hbox{with} \qquad
     u =\varphi   \ \hbox{  on  } \partial D.
 $$
 Here  $X$ is a two-dimensional Brownian motion,
  $\P_x$ is the law of $X$ starting from $x$ at time $0$,  $\E_x$ denotes the mathematical expectation under  the probability measure  $\P_x$, 
  and $\tau_D:=\inf\{t>0: X_t \notin D\}$ is the first exit time from $D$ by $X$.
In this paper, we use $:=$ as a way of definition.
Since then, there are many subsequent works using Markov processes to solve partial differential and integro-differential equations;
see, e.g.,  \cite{CZ1, CZh1, CZ, MS, Su} and the references therein.
For instance, let
\begin{equation}\label{e:1.1}
\LL^0 =  {\rm div} ( A(x) \nabla) + b(x)  \cdot \nabla,
\end{equation}
where $A(x)=(a_{ij}(x))_{1\leq i, j\leq d}$ is $d\times d$-matrix-valued measurable function on $\R^d$ that is uniformly elliptic and bounded,  that is, there is a constant $\lambda \geq 1$ so that
\begin{equation}\label{e:elliptic}
\lambda^{-1} |\xi |^2\leq \sum_{i,j=1}^d a_{ij}(x) \xi_i \xi_j \leq \lambda |\xi|^2
\quad \hbox{for every } x, \xi \in \R^d,
\end{equation}
and
$ b (x)=(b_1(x), \dots, b_d (x))$  is an $\R^d$-valued function defined on $\R^d$ so that $|b|^2$ is in the Kato class $\K_d$ (see
Definition \ref{D:2.1} below).   It is shown in \cite{CZh1} via a semi-Dirichlet form construction
that there is  a  conservative diffusion  process $X^0$ on $\R^d$
whose $L^2$-generator is $\LL^0$.
Let $D$ be a bounded domain in $\R^d$ with $d\geq 1$ which is regular with respect to the Laplace operator
(see Definition \ref{D:3.3} below),
and $ \kappa\leq 0$ be a function in $\K_d$.
 It is shown in \cite[Theorem 1.1]{CZ1} that for any continuous function $\varphi$ on $\partial D$,
 $$u(x):=\E_x [ \varphi (X^\kappa_{\tau_D})], \quad x\in D,
 $$
 is the unique weak solution to
 $  (\LL^0+\kappa)  u =0 $ in $D$ that is continuous on $\overline D$  with $u=\varphi$ on $\partial D$,
where  $X^\kappa$ is the subprocess of $X^0$ killed at rate $\kappa$. Thus the above solution $u$
can also be expressed as
$$
u(x)=\E_x \left[ \exp \left(\int_0^{\tau_D} \kappa(X^0_s) ds \right) \varphi (X^0_{\tau_D}); \tau_D<\infty \right], \quad   x\in D.
$$

When $\sL$ is an integro-differential operator, the Markov process $X$ associated with it, if there is one,
has discontinuous sample paths. Unlike the diffusion processes, the Markov process $X$
may leave a domain $D$ by  jump across the boundary $\partial D$. So for an
integro-differential operator  $\sL$, its Dirichlet problem should be an exterior value problem; that is,
 $$
   \sL  u =0  \  \hbox{  in } D  \qquad \hbox{with} \qquad
     u =\varphi   \ \hbox{  on  }  D^c.
 $$
It is natural to guess that  the solution to the above Dirichlet problem can be given probabilistically by
\begin{equation} \label{e:1.3}
u(x)=\E_x [ \varphi (X_{\tau_D}); \tau_D <\infty], \quad x\in D.
\end{equation}
The central and the hard part,   for  both differential   and   integro-differential operators,
    is to show that the function $u$ defined  by  \eqref{e:1.3} with $\varphi \in C_b (\partial D)$ or $\varphi \in C_b (D^c)$
    does give a weak solution for the Dirichlet problem in the analytic sense and such a weak solution is unique.

When the bilinear form generated by $\sL$ is a symmetric  regular Dirichlet form or non-symmetric  regular  Dirichlet form
satisfying the sector condition and $\varphi$ is a function in the domain of the Dirichlet form, it follows directly from the Dirichlet form
theory, see \cite{CF, FOT, MR}, that the Dirichlet problem for $\sL$ can be solved by a suitable projection operator
and it admits the probabilistic representation
\eqref{e:1.3}.  For a general bounded $\varphi$ on $D^c$, it is established in \cite{CH} that the function $u$ given by  \eqref{e:1.3} is an $\sL$-harmonic function in $D$ in the analytic sense   when $\sL$ is  the generator of  a symmetric  regular Dirichlet form.
This result has been extended in \cite{MZZ} to  non-symmetric operator $\sL$ that  is the generator of a non-symmetric sectorial
 regular Dirichlet forms.

In this paper, we study the Dirichlet problem on bounded regular domains $D\subset \R^d$ for the following integro-differential operator
\begin{equation} \label{e:1.4}
\sL u (x):= \sL^0 u(x) + \int_{\R^d} (u(y)-u(x)) J(x, dy),
\end{equation}
where $\sL^0$ is the operator in \eqref{e:1.1}, and   $J(x, dy)$ is a finite measure on $\R^d$ for each $x\in \R^d$
with $\kappa (x):=J(x, \R^d)\in  \K_d$.
Without loss of generality, we may  and do assume $J(x, \{x\})=0$.
Note that the   measure $J(x, dy)$ is quite general and may have atoms so it may charge on polar sets of Brownian motion or of the diffusion
process associated with the elliptic operator $\sL^0$  on $\R^d$.
In general, the generator $\sL$ is non-symmetric and its bilinear form does not satisfy the sector condition
so Dirichlet form approach is not directly applicable here.

\smallskip

 In this paper, for a Borel set $A\subset \R^d$,  $C_b(A)$ denotes the space of bounded continuous functions on $A$.
 We   denote by   $C_\infty (\R^d)$  the space of continuous functions on $\R^d$ that vanishes at infinity
equipped with the Banach norm $\| f\|_\infty := \sup_{x\in \R^d} |f(x)| $.
 Recall that a Markov process $X=\{X_t, t\geq 0; \P_x, x\in \R^d\}$ is said to be a Feller process if its   transition semigroup $\{P_t; t\geq 0\}$
is a strongly continuous semigroup in the Banach space $(C_\infty (\R^d), \| \cdot \|_\infty)$,
 and it is said to have strong Feller property if  for each $t>0$, $P_t$  maps bounded functions on $\R^d$
to functions in $C_\infty (\R^d)$. We are concerned with the following two questions in this paper.

\begin{enumerate}
\item[(i)]  Is there a Markov process on $\R^d$ associated with the non-local operator $\sL$ and in which sense? How to characterize
such a Markov process if there is one?

\item[(ii)]  How to solve the Dirichlet problem for $\sL$ on bounded domains $D\subset \R^d$? Is the solution unique and in which sense?
Note that since the jump measure $J(x, dy)$ in \eqref{e:1.4} may have atoms, the solution, if exists, needs to be uniquely defined everywhere in $D$.
\end{enumerate}

Below is the summary of the main results and contributions of this paper.

\begin{enumerate}
\item[(i)] The first main result of this paper is Theorem \ref{T:2.5}, in which we show there is a unique Feller process
$X=\{X_t, t\geq 0; \P_x, x\in \R^d\}$ having strong Feller property
associated with generator $\sL$ of \eqref{e:1.4}.
 While the desired Markov process
$X$ can be obtained from the diffusion process $X^0$ associated
with $\sL^0$ through a redistributing or piecing together procedure, the main part of the work is to establish the uniqueness of
the Markov process associated with $\sL$ and to show that the  process $X$ is a Feller process having strong Feller
property. We point out that due to the singular nature of the drift $b$ and the jump measure $J(x, dy)$,
it is difficult to employ the martingale problem approach and the SDE approach.
We formulate and established the existence and uniqueness problem in terms
of the generator and resolvents.

\item[(ii)] The second main result of this paper is Theorem \ref{T:3.4},  which says that for any bounded regular domain
$D\subset \R^d$ and any bounded continuous function $\varphi $ on $D^c$,  $u(x):=\E_x [\varphi (X_{\tau_D})]$
is the unique weak solution to $\sL u=0$ that is continuous on $\R^d$ with $u =\varphi $ on $D^c$.
The crux work is to show that $u(x):=\E_x [\varphi (X_{\tau_D})]$ is  weak differentiability in $D$ and  continuous on $\R^d$, and
it uniquely solves the Dirichlet problem in the analytic sense, which is  non-trivial.  As mentioned earlier, the bilinear form associated with $\sL$
typically does not satisfy the sector condition. Moreover, $J(x, dy)$ may have atoms.
So we can not employ the Dirichlet form method to establish the regularity of $u$.
Our method extends easily to solve the Poisson-Dirichlet problem for Schr\"odinger operator
$$
  ( \sL  +q) u =f  \  \hbox{  in } D  \qquad \hbox{with} \qquad
     u =\varphi   \ \hbox{  on  }  D^c,
 $$
 where $q$ and $f$ are function in $\K_d$ with $q^+:=\max\{q, 0\}$ satisfying certain smallness condition (see \cite{CZ}).
 But for the transparency of the exposition, in this paper we confine ourself to the case that $q=f=0$ in $D$.
\end{enumerate}

While there are many works on using probabilistic method solving Dirichlet problem for differential operators and some works
for purely non-local operators, there is very limited work
for integro-differential operators that have diffusive parts. In  \cite{Su},   it is shown, using the two-sided heat kernel estimates from \cite{CH} and 
the  semi-Dirichlet form theory,
that for a bounded Lipschitz domain  $D$ of $\R^d$,
  the following Dirichlet (complement value) problem
$$
   \begin{cases}
    ( \Delta+ b\cdot \nabla  + a \Delta^{{\alpha}/2}  +c) u   =f  \quad  &\hbox{in}  D ,  \\
     u =\varphi    \quad &\hbox{on }   D^c , 
   \end{cases}
   $$
 has a unique bounded weak solution $u$ that is continuous in $D$ and it has a probabilistic representation,
 where $a>0$ is a constant, $|b|^2$, $c$ and $f$ satisfy certain $L^p$-integrability condition   which are stronger than the Kato class condition $\K_d$
 with   small $L^p$-norm condition on $c^+:=\max\{c, 0\}$, 
 and $\varphi$ is a bounded function on $D^c$.  In the above, $\Delta^{{\alpha}/2}$ is the fractional Laplacian on $\R^d$ with $\alpha \in (0, 2)$,
 which is the generator of an isotropic $\alpha$-stable process on $\R^d$. The fractional Laplacian
   can be expressed as
 $$
 \Delta^{\alpha/2} \psi (x) = \lim_{\eps \to 0+} \int_{\{ y\in \R^d: |y-x|\geq  \eps\}}  (\psi (y)-\psi (x)) \frac{c_{d, \alpha}}{|y-x|^{d+\alpha}} dy
 $$
 for $C^2$-smooth functions with bounded derivatives, where $c_{d, \alpha}>0$ is an explicit constant that depends on $d$ and $\alpha$.
  Thus $\Delta^{\alpha/2} \psi (x)$ has a jump kernel
 $J(x, dy) = \frac{c_{d, \alpha}}{|y-x|^{d+\alpha}} dy$, which is an infinite measure on $\R^d$ of a particular form and 
 the measure $J(x, dy) dx$ is    symmetric   on $\R^d\times \R^d$.   
 Whether   the weak solution $u$ is continuous up to the boundary $\partial D$
 for  $\varphi\in C_b(D^c)$, however,  is not studied in \cite{Su}.

\smallskip

In this paper, we focus   on the  Dirichlet problem for integro-differential operators of \eqref{e:1.4}
and its   probabilistic representation, where $J(x, dy)$
 in \eqref{e:1.4} is a finite (arbitrary) measure  on $\R^d$ for each $x \in \R^d$.
  In a subsequent work, we will study the Dirichlet problem for $\sL$ of \eqref{e:1.4} where the kernel $J(x, dy)$ can   be an infinite measure
 which can also be   degenerate  but $J(x, dy) dx$ is symmetric in $(x, y)$.
  We refer the reader  to \cite{CZ2, Hs, Pa} and the references therein   for solutions to  differential operators  with other boundary condition and
their probabilistic representations, and to \cite{AR, FKV} and the references therein for analytic approaches to Dirichlet problem
for purely  non-local operators.

\smallskip

 The rest of the paper is organized as follows.  In Section \ref{S:2}, we give definitions of Kato class and weak solutions, and present the proof of the
 first main result of the paper  after some preparations.
 The proof of the second main result of this paper is given in Section \ref{S:3}.

\section{Diffusion with jumps}\label{S:2}

Let $d\geq 1$ and $B=\{B_t, t\geq 0; \P_x, x\in \R^d\}$ be Brownian motion on $\R^d$, where $\P_x$ is the law of $B$ starting from $x$ at $t=0$.

\begin{defn}\label{D:2.1}  \rm
A real-valued function $f$ on $\R^d$
is said to be in Kato class $\K_d $
if
\begin{equation}\label{e:2.1}
\lim_{t\to 0} \sup_{x\in \R^d}  \E_x \int_0^t |f(B_s)| ds =0 ,
\end{equation}
where   $\E_x$ denotes the mathematical expectation under  the probability measure  $\P_x$.
\end{defn}

It is known   that $f\in \K_d$ if and only if
\begin{equation}\label{e:2.2}
\begin{cases}
\lim_{r\to 0} \sup_{x\in \R} \int_{B(x, 1)}  |f(y)|  dy <\infty \quad &\hbox{when } d=1,  \smallskip \cr
  \lim_{r\to 0} \sup_{x\in \R^2} \int_{B(x, r)}  |f(y)| \,  \log (1/|x-y|) dy =0
  \quad &\hbox{when } d=2,   \smallskip \cr
\lim_{r\to 0} \sup_{x\in \R^d} \int_{B(x, r)} \frac{ |f(y)|}{|x-y|^{d-2}} dy =0
\quad & \hbox{when } d\geq 3 .
\end{cases}
\end{equation}
Here $B(x, r)$ denotes the ball in $\R^d$ centered at $x \in \R^d$ with radius $r > 0$.
It is also known that
$$
L^{\infty} (\R^d) + L^p(\R^d)   \subset \K_d \quad \hbox{for every } p>d/2.
$$
See, e.g.,  \cite{CZ} for above assertions. Observe that  \eqref{e:2.2} implies that $\K_d \subset L^1_{\rm loc} (\R^d)$.
For simplicity, in the following we refer the least upper bound of the convergence rate
in \eqref{e:2.1}, or equivalently, in \eqref{e:2.2}, as the Kato norm of $f$.

\smallskip

Let
\begin{equation}\label{e:2.3}
\LL^0 =\frac{1}{2} {\rm div} ( A(x) \nabla) + b(x)  \cdot \nabla,
\end{equation}
where $A(x)=(a_{ij}(x))_{1\leq i, j\leq d}$ is $d\times d$-matrix-valued measurable function on $\R^d$
that is uniformly elliptic and bounded (that is,  it satisfies condition \eqref{e:elliptic}),
  and
$ b (x)=(b_1(x), \dots, b_d (x))$  is an $\R^d$-valued function defined on $\R^d$ so that $|b|^2\in \K_d$.
Here we put $\frac12$ in \eqref{e:2.3} just for convenience as the infinitesimal generator of Brownian motion is $\frac12 \Delta$.
In \cite{CZh1}, it is shown via a semi-Dirichlet form construction that there is a continuous conservative Feller  process
$X^0$ on $\R^d$ having strong Feller property having $\LL^0$ as its $L^2$-generator.

\smallskip

Suppose that for each $x\in \R^d$, $J(x, dy)$ is a finite measure on $\R^n$
with $J(x, \{x\})=0$ and $\kappa (x):=J(x, \R^d)$ being   in the Kato class $\K_d$.
We assume $\kappa$ is non-trivial and consider the following operator
\bea  \label{e:2.4}
\LL u(x) =\LL^0 u (x) +  \int_{\R^d}  ( u(y) -u(x) ) J(x, dy),
  \eea
on $\R^d$. Define $\nu (x, dy)= J(x, dy)/\kappa (x)$.
Clearly, $\nu (x, \{ x\})=0$ for every $x\in \R^d$.
We can rewrite \eqref{e:2.4} as
\bea
\LL u(x) =\LL^0 u (x) +\kappa (x)  \int_{\R^d}  ( u(y) -u(x) ) \nu (x, dy).
  \label{e:2.5}
  \eea
The main result of this section is Theorem \ref{T:2.5}, which asserts that there is a unique Feller process $X$ having $\sL$ as its generator.
We will further show that $X$  has strong Feller property and it can be constructed from the diffusion process $X^0$ through a redistributing or
piecing-together procedure.  To present the theorem precisely, we need some notations and preparations.

\smallskip

For an open subset $D\subset \R^d$,
 denote  by $C^1_c(D)$ the space of continuously differentiable functions with compact support in $D$. Define
 $$
W^{1,2}(D)=\{f\in L^2(D; dx):\nabla f\in L^2(D; dx)\},
$$
where the gradient $\nabla f$ is understood in the distributional sense. A function $f$ is said to be
locally in $W^{1,2}(D)$, denoted as $f\in W^{1,2}_{\loc }(D)$,
 if $f\in W^{1,2}(D_1)$ for every relatively compact subdomain $D_1$ of $D$.
The closure of 
$C^1_c(D)$  with respect to the Sobolev norm
\begin{equation}\label{e:2.6}
\| f\|_{1,2}:= \left(\| f \|_{L^2(D; dx)} + \| \nabla f \|_{L^2(D; dx)}\right)^{1/2}
\end{equation}
is denoted as  $W^{1,2}_0(D)$. Intuitively speaking, $W^{1,2}_0(D)$ consists of those functions in $W^{1,2}(D)$
that vanish on $\partial D$.
 One of the important properties of a function $q$  in the Kato class $\K_d$ is that for any
$\eps > 0$, there exists a constant $A_\eps > 0$  such that (cf. \cite{Kato}) for any $u\in W^{1,2}(\R^d)$,
$$
\int_{\R^d} u(x)^2 |q(x)| dx \leq \eps \int_{\R^d} |\nabla u (x)|^2 dx+ A_\eps \int_{\R^d} u(x)^2 dx.
$$

\begin{defn}\label{D:2.2} \rm
Let $D$ be a connected open subset of $\R^d$ and $q$ a real-valued function on $D$ so that $1_D q \in \K_d$.
Let $f$ be a bounded function on $D$.
\begin{description}
\item{(i)}
A function $u$ on $D$ is said to be a weak solution of $(\LL^0 + q) u=f$ in $D$ if
$u\in W^{1,2}_{\loc }(D)$ and
\begin{eqnarray*}
&& -\frac{1}{2}\sum_{i,j=1}^d\int_D a_{ij} (x) \frac{\partial u (x) }{\partial x_i}\frac{\partial \psi (x) }{\partial x_j}dx
+ \sum_{i=1}^d \int_D b_i(x)  \frac{\partial u (x) }{\partial x_i} \psi (x) dx
 + \int_D q(x) u (x) \psi (x) dx \\
 &=&  \int_D f(x) \psi (x) dx
\end{eqnarray*}
for every $\psi \in C_c^1(D)$.
\end{description}

\item{(ii)} A   function $u$ on $\R^d$ is said to be a weak solution of $(\LL + q) u=f$ in $D$ if
$u\in W^{1,2}_{\loc }(D) \cap C_b (D)$ and
\begin{eqnarray*}
&&    -\frac{1}{2}\sum_{i,j=1}^d\int_D a_{ij} (x) \frac{\partial u (x) }{\partial x_i}\frac{\partial \psi (x) }{\partial x_j}dx
+ \sum_{i=1}^d \int_D b_i(x)  \frac{\partial u (x) }{\partial x_i} \psi (x) dx
 + \int_D q(x) u (x) \psi (x) dx\\
&& \quad  + \int_D \left( \int_{\R^d}  ( u(y) -u(x) ) J(x, dy)\right) \psi (x) dx= \int_D f(x) \psi (x) dx.
\end{eqnarray*}
for every $\psi \in C_c^1(D)$.

\item{(iii)} A function $u$ is said to be $(\sL^0+q)$-harmonic (respectively, $(\LL + q)$-harmonic) in $D$ if
it is a weak solution to $(\sL^0+q) u=0$ (respectively, $(\LL + q) u=0$) in $D$.
\end{defn}

 Note that in (ii) of the above definition for the weak solution of $(\LL^0 + q) u=f$ in $D$,
we require $u$ to be bounded and continuous on $D$ because $J(x, dy)$ can possibly charge on singletons.

\smallskip

 Let $X^\kappa$ be the subprocess of $X^0$ killed at (state-dependent) rate $\kappa (x)$;
that is,
$$
\E_x [ f(X^{\kappa}_t)] = \E_x  [ e_\kappa (t) f(X_t^0)], \quad t\geq 0, \  x\in \R^d,
$$
where $e_\kappa (t):=\exp\left(- \int_0^t \kappa (X^0_s) ds \right)$.
Since $X^0$ is a Feller process on $\R^d$ having strong Feller property and $\kappa \in \K_d$,
it follows from \cite[Theorem 2]{Chu} that $X^\kappa $ is a Feller process on $\R^d$ having strong Feller property.
 It is established in \cite{CZh1} that $X^\kappa$ has $L^2$-infinitesimal  generator $\LL^0-\kappa$.
Denote by $\zeta^\kappa$ the lifetime of $X^\kappa$.
For $\alpha\geq 0$, the $\alpha$-resolvent  operator $G^\kappa_\alpha$ of $X^\kappa$ is given by
$$
G^\kappa_\alpha f(x):= \E_x \int_0^\infty e^{-\alpha t} f(X^\kappa_t) dt, \quad x\in \R^d,
$$
for any $f\geq 0$ on $\R^d$. The $\alpha$-resolvent $G^0_\alpha$ for $X^0$ is defined in an analogous  way.
The following result is proved in \cite[Proposition 2.2]{CZh2}.

\begin{prop}\label{P:2.3}
For every $\varphi \geq 0$ and $\alpha \geq 0$,
$$
\E_x \left[ e^{-\alpha \zeta^\kappa} \varphi (X^\kappa_{\zeta^\kappa -} ); \zeta^\kappa <\infty \right]
=G^\kappa_\alpha (\kappa \varphi )(x), \quad x\in \R^d.
$$
\end{prop}

\smallskip

We recall from \cite{CZh1} another construction of $X^0$ through the Girsanov transform of the  symmetric symmetric diffusion process
 $\bar X^0=\{\bar X^0, t\geq 0; \bar \P^0_x, x\in \R^d\}$ on $\R^d$ having $\frac12 \nabla (A(x) \nabla)$ as its infinitesimal generator.
It is well known \cite{An} that $\bar X^0$ has a jointly continuous transition density function $\bar p_0 (t, x, y)$ on $(0, \infty) \times \R^d \times \R^d$
that admits the following two-sided Aronson estimate:
\begin{equation}\label{e:2.7}
c_1 t^{-d/2} e^{-c_2|x-y|^2/t} \leq \bar p_0(t, x, y)\leq c_3 t^{-d/2} e^{-c_4 |x-y|^2/t}
\quad \hbox{for all } t>0 \hbox{ and } x, y\in\R^d.
\end{equation}
The symmetric diffusion $\bar X^0$ is not a semimartingale in general. It is a Dirichlet process whose martingale part is $\int_0^t \sigma (\bar X^0_s) dB_s$, where $\sigma (x)$ is the square root matrix of $A(x)$ and $B$ is a $d$-dimensional Brownian motion; see \cite[(2.12)]{CZh1}.
Define
\begin{equation}\label{e:2.8}
M_t = \exp \left( \int_0^t (\sigma^{-1}b)(\bar X^0_s)dB_s - \frac12 \int_0^t |(\sigma^{-1}b)(\bar X^0_s)|^2 ds \right), \quad t\geq 0.
\end{equation}
It is shown in \cite[Theorem 3.1]{CZh1} that $M:=\{M_t, t\geq 0\}$  is a square integrable martingale with respect to the augmented filtration generated by $\bar X^0$ with
\begin{equation}\label{e:2.9}
 C_0:=\sup_{x\in \R^d} \bar \E^0_x \left[ M^2_{t_0}\right] <\infty
\end{equation}
for some $t_0>0$. In fact, by the Markov property of $\bar X^0$, the above property holds for every $t_0>0$.
Furthermore, it is established in \cite[Theorem 3.2]{CZh1} that the diffusion process
$X^0=\{ X^0, t\geq 0;  \P^0_x, x\in \R^d\}$
 with $L^2$-infinitesimal generator $\LL^0=\frac12 \nabla \cdot (A(x) \nabla)+b \cdot \nabla$
has the same distribution as
the Girsanov transformed process obtained from $\bar X^0$ through the exponential martingale $M$.
That is,
\begin{equation} \label{e:2.10}
\E_x \left[ f(X^0_t)\right] = \bar \E_x [ M_t f(\bar X^0_t) ], \quad x\in \R^d, \ t\geq 0,
\end{equation}
for every $f\geq 0$ on $\R^d$.
Recall that for an open set $U\subset \R^d$, we  use the   notation $\tau_U$
to denote the first exit time from $U$ by $X^0$  as well as by other processes, say, $\bar X^0$,
whose meaning should be clear from the context.

\begin{lemma}\label{L:2.4}  For every  $q\in \K_d$, there exists a constant $\alpha_0>0$ that depends only on the Kato norm of $q$ so that
for every $\alpha \geq \alpha_0$,  both $G^0_\alpha q $ and $G^\kappa_\alpha q $ are in $ C_\infty (\R^d)$.
 \end{lemma}

\pf    Let   $q\in \K_d$.  By  \eqref{e:2.9}-\eqref{e:2.10}, Cauchy-Schwarz inequality and \eqref{e:2.7},
\begin{eqnarray}\label{e:2.11}
\lim_{t \to 0} \sup_{x\in \R^d} \E_x  \int_0^t  | q   (X^\kappa_s)  | ds
&\leq & \lim_{t \to 0} \sup_{x\in \R^d} \E_x^0 \int_0^t  | q (X^0_s)|  ds  \nonumber \\
&=&  \lim_{t \to 0} \sup_{x\in \R^d}  \bar \E_x^0  \left[M_t \int_0^t  | q  (\bar X^0_t) |  ds \right] \nonumber \\
&\leq & c \lim_{t \to 0} \left( \sup_{x\in \R^d}  \bar \E_x^0  \left[ \int_0^t  | q (B_t) | ds \right]  \right)^{1/2}=0.
\end{eqnarray}
Hence there is some  $\alpha_0>0$ that depends only on the Kato norm of $q$ so that
  $G^0_{\alpha_0} |q|$ and hence $G^\kappa_{\alpha_0} |q| $ is bounded on $\R^d$.
 Denote by $\{P^0_t, t\geq 0\}$ and $\{P^\kappa_t, t\geq 0\}$ the transition semigroups of $X^0$ and $X^\kappa$, respectively.
Since $X$ and $X^\kappa$ are  Feller processes having strong Feller property,
$P^0_t (G^0_\alpha q)\in C_\infty (\R^d)$ and
$P^\kappa_t (G^\kappa_\alpha q)\in C_\infty (\R^d)$ for every $t>0$ and $\alpha \geq \alpha_0$.
Note that
\begin{equation}\label{e:2.18}
 G_\alpha^\kappa q(x)  = \int_0^t e^{-\alpha s}  P^\kappa_s q (x)ds  + e^{-\alpha t} P^\kappa_t (G_\alpha^\kappa q)(x)
\end{equation}
  This combined with \eqref{e:2.11} implies that $G^\kappa_\alpha  q \in C_\infty (\R^d)$ for every $\alpha\geq \alpha_0$.
The same reasoning shows that  $G^0_\alpha  q \in C_\infty (\R^d)$ for every $\alpha\geq \alpha_0$.
\qed

 We now present  the main result of this section.

\begin{thm}\label{T:2.5}
\begin{enumerate}
\item [\rm (i)]  There is a conservative Hunt process $X=\{X_t, t\geq 0; \P_x, x\in \R^d\}$ having  $\LL$ as its generator  in the following sense.
There is some $\alpha_0>0$ that depends only on the ellipticity constant $\lambda$ of $A(x)=(a_{ij}(x))$ and the Kato class norm of $|b|^2 + \kappa $ so that for every $\alpha >\alpha_0$ and for any
bounded function $f$ on $\R^d$,
$G_\alpha f(x):= \E_x \left[ \int_0^\infty e^{-\alpha t} f(X_t) dt\right]$ is in $W_{\loc }^{1,2}(\R^d) {\cap C_b(\R^d)}$
and $(\alpha - \LL) G_\alpha f= f $ in the weak  sense.

\item [\rm (ii)] Suppose that $Y$  is another Hunt process  on $\R^d$  having generator $\LL$  (but without assuming a priori the conservativeness).
Then $Y$ has the same distribution as $X$. Consequently,   $Y$ has infinite lifetime.

\item[\rm (iii)]  The conservative Hunt process $X=\{X_t, t\geq 0; \P_x, x\in \R^d\}$ having generator  $\LL$ is a Feller process on $\R^d$
having strong Feller property.
  \end{enumerate}
\end{thm}

\pf  (i)
 We  can  construct a strong Markov process $X$ from the diffusion process $X^0$ through a re-distributing or piecing-together procedure as follows.
Recall that  $X^\kappa$ is  the subprocess of $X^0$ killed at (state-dependent) rate $\kappa (x)$.
  Given any initial starting point $x_0\in \R^d$, we run a copy of $X^\kappa$ starting from $x_0$.
 At its death place, say $x_1$,
re-start an independent copy of $X^\kappa$ starting from a  random point $y_1$ selected according to the probability measure
$\nu (x_1, dy)$. At the next death place, say, $x_2$, start another independent copy of $X^\kappa$ starting from a random point $y_2$ selected
according to the probability measure $\nu (x_2, dy)$.
Continuing  this procedure results a Hunt process $X=\{X_t, t\geq 0; \P_x, x\in \R^d\}$; see \cite{INW}.
We next show that $X$ is conservative.
Let
$$
\tau:=\inf\{t>0: X_t \not= X_{t-}\},
$$
which is the lifetime of $X^\kappa$. We know by \cite[Theorem 5.12]{CZh1} by taking $D=\R^d$ there
that
\begin{equation}\label{e:2.13}
\lim_{\alpha \to \infty} \sup_{x\in \R^d} \E_x \left[ e^{-\alpha \tau} \right] =0.
\end{equation}
So there is $\alpha >0$ so that
\begin{equation}\label{e:2.14}
\sup_{x\in \R^d} \E_x \left[ e^{-\alpha \tau} \right] \leq 1/2.
\end{equation}
Define $\tau_0=0$, $\tau_1=\tau$, $\tau_{k+1}= \tau_k + \tau\circ \theta_{\tau_k}$ for $k\geq 1$;
that is, $\tau_k$ is the $k$th re-distributing time. By the strong Markov property of $X$,
$$
\sup_{x\in \R^d} \E_x \left[ e^{-\alpha \tau_k} \right]
= \sup_{x\in \R^d} \E_x \left[ \prod_{j=1}^k e^{-\alpha (\tau_j-\tau_{j-1})} \right]
\leq \left( \sup_{x\in \R^d} \E_x \left[ e^{-\alpha \tau} \right]\right)^k
\leq 2^{-k}.
$$
Hence
$$
\E_x \left[ e^{-\alpha \lim_{k\to \infty} \tau_k}\right]
= \lim_{k\to \infty} \E_x \left[ e^{-\alpha \tau_k} \right]=0.
$$
It follows that $\P_x (\lim_{k\to \infty} \tau_k=\infty)=1$.
This implies that the process $X$ can only have finite many re-distributing (or re-generations) during any finite time interval
and so $X$ has infinite lifetime.

Denote by $G_\alpha$ and $G^\kappa_\alpha$ the resolvent of $X$ and $X^\kappa$, respectively.
Applying the strong Markov property at the first jumping time $\tau$ of $X$,
we have by the construction of $X$ that for every $\alpha \geq 0$ and any bounded function $f\geq 0$ on $\R^d$,
\begin{eqnarray}
G_\alpha f(x)&:=& \E_x \left[\int_0^\infty e^{-\alpha s} f(X_s) ds \right] \nonumber \\
&=& G^\kappa_\alpha f(x) + \E_x \left[ e^{-\alpha \tau} G_\alpha f(X_{\tau})\right] \nonumber \\
&=&  G^\kappa_\alpha f(x) + \E_x \left[ e^{-\alpha \tau} \int_{\R^d}
G_\alpha  f(y) \nu (X_{\tau-}, dy)\right] \nonumber \\
&=&  G^\kappa_\alpha f(x) +  G^\kappa_\alpha \left(\kappa (\cdot) \int_{\R^d}
G_\alpha f (y)  \nu (\cdot, dy) \right) (x) \nonumber \\
&=&  G^\kappa_\alpha f(x) +  G^\kappa_\alpha \left( \int_{\R^d}
G_\alpha f (y)  J (\cdot, dy) \right) (x) ,  \label{e:2.15}
\end{eqnarray}
where the second to the last equality is due to Proposition \ref{P:2.3}.
For any relatively compact subset $D\subset \R^d$,
denote by $\tau_D$ the first exit time from $D$ by $X^\kappa$.
For $\alpha \geq 0$ and $g\geq 0$ on $D$, define
$$
G^{\kappa, D}_\alpha g(x):= \E_x \left[ \int_0^{\tau_D} e^{-\alpha t} g(X^\kappa_t) dt \right].
\quad x\in D,
$$
 Let $g(x):= f(x) + \int_{\R^d} G_\alpha f(y) J(x, dy)$, whose absolute value is  bounded by $\| f \|_\infty (1+ \kappa (x)  /\alpha)$ on $\R^d$.
We have by \eqref{e:2.15} and the strong Markov property of $X^\kappa$ at $\tau_D$ that
\begin{equation}\label{e:2.16}
G_\alpha f(x) = G_\alpha^{\kappa } g(x)   = G_\alpha^{\kappa, D} g(x)
       + \E_x \left[ e^{-\alpha \tau_D} G^\kappa_\alpha g(X^\kappa_{\tau_D}) \right],
\quad x\in D.
\end{equation}
Here    $G_\alpha^{\kappa, D}$ is the $\alpha$-resolvent of the subprocess of $X^\kappa$ killed upon leaving $D$.

Let $\alpha_0>0$ be the constant in Lemma \ref{L:2.4} for $q=\kappa + 1$.
Note that $f\in L^2(D)$. By \cite[Theorem 5.2]{CZh1}, there is some $\alpha_1\geq \max\{\alpha_0, 1\}$
that depends only on the ellipticity bounds of $A(x)$ and the Kato class norm of $|b|^2 + \kappa $ so that for every $\alpha > \alpha_1$,
$G^{\kappa, D}_\alpha f \in W^{1,2}_0(D)$ and $(\alpha + \kappa -\LL^0)G_\alpha^{\kappa, D}f=f$ in the weak (distributional) sense.
On the other hand, as  $|g (x)| \leq  \| f \|_\infty (1+   \kappa (x) /\alpha)$ on $\R^d$,    $  G_\alpha^\kappa g \in C_\infty (\R^d)$
for $\alpha >\alpha_1$ by Lemma \ref{L:2.4}.
Now by \cite[Theorem 5.11]{CZh1} with $q=-(\kappa +\alpha)$ there that
$$ v(x):=\E_x \left[ e^{-\alpha \tau_D} G^\kappa_\alpha  g(X^\kappa_{\tau_D}) \right]
 = \E_x^0 \left[ e_{-\kappa - \alpha} (\tau_D) G^\kappa_\alpha  g(X^0_{\tau_D}) \right] , \quad x\in D,
$$
is in $W^{1,2}_{\loc}(D) \cap C(D)$
 and $(\LL^0-\kappa - \alpha)$-harmonic in $D$.
We conclude from \eqref{e:2.16} that for every $\alpha >\alpha_1$,
$G_\alpha f \in W^{1,2}_{\loc} (\R^d) \cap C_\infty (\R^d)$
and that
$$
(\alpha+ \kappa-\LL^0) G_\alpha f = g = f(x) + \int_{\R^d} G_\alpha f(y) J(x, dy)
$$
in the weak sense. The latter is equivalent to
$(\alpha-\LL) G_\alpha f = f(x)$ in the weak sense.
This establishes that $X$ is a conservative Hunt process associated with $\sL$.

\smallskip

(ii) We next show   the uniqueness of the Hunt process with generator $\LL$.
Suppose that $Y=\{Y_t, t\geq 0; \P_x, x\in \R^d\}$ is a Hunt process so that
there is some $\alpha_2 >0$  such  that for every $\alpha > \alpha_2$ and for bounded function $f$ on $\R^d$,
$G^Y_\alpha f(x):= \E_x \left[ \int_0^\infty e^{-\alpha t} f(Y_t) dt\right]$ is in $W_{\loc }^{1,2}(\R^d)$
and $(\alpha - \LL) G^Y_\alpha f= f $ in the weak  sense.
Note that we do not assume a priori that $Y$ has infinite lifetime.
Without loss of generality, we may and do assume that $\alpha_2 \geq \alpha_1$.
Then for every $\alpha \geq \alpha_1$ and for bounded function $f$ on $\R^d$,
 $u=G^Y_\alpha f -G_\alpha f  \in C_b(\R^d)$
is a weak solution to $(\alpha -\LL)u=0$ on $\R^d$, or equivalently,
$$
(\alpha + \kappa -\LL^0) u =  \int_{\R^d} u(y) J(x, dy)=:g(x)
\quad \hbox{on } \R^d.
$$
Clearly,  $u$ is bounded and
$ |g(x)|=\left| \int_{\R^d}  u (y) J(x , dy) \right| \leq \| u\|_\infty \kappa (x). $
 For an integer $k\geq 1$, let $B_k:=B(0, k)$ and
$v_k (x):=\E_x \left[ e^{-\alpha \tau_{B_k}} u(X^\kappa_{\tau_{B_k}}); \tau_{B_k}  <\zeta^\kappa \right]$.
Here $\tau_{B_k}:=\inf\{t>0: X^\kappa_t \notin B_k\}$.
We know by
\cite[Lemmas  5.6 and 5.8]{CZh1}
 that
$v_k$ is the weak solution of $(\alpha +\kappa -\LL^0)v=0$
in $B_k$
so that
$v_k-1_{B_k} u\in W^{1, 2}_0(B_k)$. Hence $u_k:= u-v_k \in W^{1, 2}_0 (B_k)$ is a weak solution
to
$ (\alpha + \kappa -\LL^0) u_k = g$ in $B_k$.
Thus  by the proof of \cite[Theorem 3.2]{CZh2}  (from line -5 on p.307 to  (3.9) on p.308 there),
the above equation has a unique weak solution
and
$$
u(x)-\E_x \left[ e^{-\alpha \tau_{B_k}} u(X^\kappa_{\tau_{B_k}})\right]
= u_k(x) = \E_x \int_0^{\tau_{B_k}} e^{-\alpha t} g(X^\kappa_t) dt, \quad x\in B_k.
$$
Noting that $\lim_{k\to \infty} \tau_{B_k}=\zeta^\kappa$ and
$\cup_{k=1}^\infty \{\tau_{B_k} = \zeta^\kappa\}=\{ \zeta^\kappa <\infty\}$,
we have by taking $k\to \infty$ in above display
$$
u(x) = \E_x \int_0^\infty e^{-\alpha t} g(X^\kappa_t) dt=
G^\kappa_\alpha g (x), \quad x\in \R^d.
$$
By Proposition \ref{P:2.3},
\begin{equation}\label{e:2.17}
|u(x)| \leq \| u\|_{\infty} G^\kappa_\alpha \kappa (x) = \| u\|_\infty \E_x
\left[ e^{-\alpha \zeta^\kappa} \right]
\quad \hbox{for every } x\in \R^d.
\end{equation}
By \eqref{e:2.14}, there is some $\alpha_3\geq \alpha_2$ so that
$\sup_{x\in \R^d} \E_x  \left[ e^{-\alpha_3 \zeta^\kappa} \right]\leq 1/2$.
It follows from \eqref{e:2.17} that for every $\alpha > \alpha_3$,
$u(x)=G^Yf-G_\alpha f$ satisfies $\| u\|_\infty \leq \| u\|_\infty /2$
and so $u\equiv 0$. By the uniqueness of the Laplace transform, we have
$\E_x [ f(Y_t)] = \E_x [ f(X_t)]$ for every $x\in \R^d$, $t>0$ and every bounded continuous function
$f\in \R^d$. Since both $X$ and $Y$ are Hunt processes, it follows that $Y$ and $X$ have the same finite dimensional distributions, and hence have the same law as processes.

\smallskip

(iii)  We now show that the conservative Hunt process $X=\{X_t, t\geq 0; \P_x, x\in \R^d\}$ associated with $\LL$ is a Feller process on $\R^d$
having strong Feller property.  For this, it suffices to show that its resolvents $\{G_\alpha, \alpha >0\}$ is a family of strongly continuous
resolvents in $(C_\infty (\R^d), \| \cdot \|_\infty)$ and for each $\alpha >0$, $G_\alpha$ maps bounded Borel measurable functions on $\R^d$
to continuous functions on $\R^d$ vanishing at infinity.
First note that as  $X^\kappa$ is a Feller process having strong Feller property,
    $G^\kappa_\alpha f \in C_\infty (\R^d)$ for every bounded Borel measurable function $f$ on $\R^d$.
On the other hand, since  $ \int_{\R^d} G_\alpha | f (y) |  J (\cdot, dy)  \leq \| f\|_\infty \kappa (x)$,
it follows from  Lemma \ref{L:2.4} that
$G^\kappa_\alpha \left( \int_{\R^d} G_\alpha f (y)  J (\cdot, dy) \right)  \in C_\infty (\R^d)$ for every $\alpha > \alpha_1$.
  Thus we have by  \eqref{e:2.15}, $G_\alpha f \in C_\infty (\R^d)$ for every bounded Borel measurable $f$ on $\R^d$.
Now suppose $f\in C_\infty (\R^d)$, we have by  \eqref{e:2.15}  and Proposition \ref{P:2.3},
$$
  \| \alpha G_\alpha f - \alpha G_\alpha^\kappa f \|_\infty
\leq   \| G_\alpha^\kappa  (1+\kappa ) \|_\infty \, \|f \|_\infty
\leq \sup_{x\in \R^d} \E_x  \left[ \alpha^{-1} (1-e^{-\alpha \zeta^\kappa})  +   e^{-\alpha \zeta^\kappa}) \right]  \| f \|_\infty.
$$
Since $X^\kappa$ is a Feller process on $\R^d$
  we have by \eqref{e:2.13}    that
 $$
 \lim_{\alpha \to \infty } \| \alpha G_\alpha f - f \|_\infty
 \leq  \lim_{\alpha \to \infty} \left( \| \alpha G_\alpha f - \alpha G_\alpha^\kappa f \|_\infty + \|\alpha G_\alpha^\kappa f  -f \|_\infty \right)
 =0.
 $$
This proves that $X$ is a Feller process on $\R^d$ having strong Feller property.
\qed

\section{Dirichlet problem and its probabilistic representation}\label{S:3}

In this section, we study the Dirichlet problem for $\sL$ on bounded regular domains by utilizing the corresponding Feller process $X$ associated with it.
Recall that from \eqref{e:2.10} that the diffusion process
$X^0=\{ X^0, t\geq 0;  \P^0_x, x\in \R^d\}$
 associated with $\LL^0=\frac12 \nabla \cdot (A(x) \nabla)+b \cdot \nabla$
can be obtained from the symmetric diffusion
$\bar X^0=\{\bar X^0, t\geq 0; \bar \P^0_x, x\in \R^d\}$ on $\R^d$ associated with $\frac12 \nabla (A(x) \nabla)$ through a Girsanov transform.

 \begin{lemma}\label{L:3.1} There are constants $c_0>0$ and $r_0>0$ that depends only on the dimension $d$, the ellipticity constant of $A$ and the Kato norm of $|b|^2$ so that for every $x_0\in \R^d$ and $r\in (0, r_0]$,
\begin{equation}\label{e:3.1}
\E^0_x \left[ \tau_{B(x_0, r)} \right] \leq c_0 r^2 , \quad x\in B(x_0, r).
\end{equation}
\end{lemma}

\pf  It follows from \eqref{e:2.7} that for any $x_0\in \R^d$, $r>0$, and $x\in B(x_0, r)$,
\begin{eqnarray*}
\bar \P^0_x (\tau_{B(x_0, r)} >t)
&\leq &   \bar \P^0_x (\tau_{B(x, 2r)} >t)
\leq  \int_{B(x, 2r)} \bar p_0(t, x, y) dy \\
c_5&\leq & c_3 \int_0^{2r/\sqrt{t}} z^{d-1} e^{-c_4 z^2} dz =: \gamma (r/\sqrt{t}).
\end{eqnarray*}
Hence by Cauchy-Schwarz inequality and \eqref{e:2.9}-\eqref{e:2.10}, for every $x\in B(x_0, r)$ and $t\in (0, t_0]$,
\begin{eqnarray*}
\P^0_x (\tau_{B(x_0, r)} >t)
 =  \bar \E^0_x \left[ M_{t}; \tau_{B(x_0, r)} >t\right]
\leq  \sqrt{ \bar \E^0_x \left[ M^2_{t} \right] \,  \bar \P^0_x ( \tau_{B(x_0, r)} >t)}
 \leq  & C_0^{1/2} \gamma (r/\sqrt{t})^{1/2}.
\end{eqnarray*}
Since $\lim_{\eps \to 0} \gamma (\eps )=\gamma (0)=0$, there is some $\eps_0>0$ so that
${ c_5} := (C_0 \gamma (\eps_0))^{1/2}<1$. Set $r_0=\eps_0 \sqrt{t_0}$.
We have from the above that for every $x_0\in \R^d$ and $r\in (0, r_0]$,
\begin{equation}\label{e:3.2}
\P_x^0 (\tau_{B(x_0, r)} > r^2/\eps_0^2) \leq  { c_5} \quad \hbox{for every } x\in B(x_0, r).
\end{equation}
Using the Markov property of $X^0$, we deduce from \eqref{e:3.2} that for every integer $k\geq 1$,
\begin{equation}\label{e:3.3}
\P_x^0 \left(\tau_{B(x_0, r)} > kr^2/\eps_0^2\right) \leq { c_5^k} \quad \hbox{for every } x\in B(x_0, r).
\end{equation}
Consequently, we get
{
$$
\E^0_x \left[ \tau_{B(x_0, r)} \right] \leq \sum_{k=1}^\infty
\frac{ kr^2}{\eps_0^2}  \, \P_x^0 \left(\tau_{B(x_0, r)} > (k-1) r^2/\eps_0^2\right)
\leq r^2 \eps_0^{-2} \sum_{k=1}^\infty   k c_5^{k-1}
= \frac{ r^2}{\eps_0^2 (1-c_5 )^2}.
$$
}
This establishes the lemma.
\qed

Recall that $X^\kappa=\{X^\kappa_t, t\geq 0; \P^\kappa_x, x\in \R^d\}$ is the subprocess of $X^0$
killed at rate $\kappa$. Here whenever needed, for emphasis, we denote the law of $X^\kappa$ with $X^\kappa_0=x$ by $\P^\kappa_x$
and its mathematical expectation by $\E^\kappa_x$. However when there is no possibility of ambiguity,
for notational simplicity,
we often drop $\kappa$ from $\E^\kappa_x$.

\begin{lemma}\label{L:3.2}
$$
\lim_{r\to 0} \sup_{x\in \R^d} \E_x \int_0^{\tau_{B(x, r)}} \kappa (X^\kappa_s) ds =0.
$$
\end{lemma}

\pf Recall that ${ c_5}\in (0, 1)$ is the constant in the proof of Lemma \ref{L:3.1} so that \eqref{e:3.2} holds.
For every $\eps_0 >0$, fix an integer $k\geq 1$ so that $2{ c_5^k}<\eps_0^2$.   By \eqref{e:3.3},
\begin{eqnarray}
&& \sup_{x\in \R^d} \E_x \int_0^{\tau_{B(x, r)}} \kappa (X^\kappa_s) ds \nonumber \\
&\leq& \sup_{x\in \R^d} \E_x \left[ \int_0^{\tau_{B(x, r)}} \kappa (X^\kappa_s) ds;
\tau_{B(x, r)} \leq kr^2/\eps_0\right]  +
 \sup_{x\in \R^d} \E_x \left[ \int_0^{\tau_{B(x, r)}} \kappa (X^\kappa_s) ds;
\tau_{B(x, r)} > kr^2/\eps_0\right]\nonumber \\
&\leq & \sup_{x\in \R^d} \E_x \left[ \int_0^{kr^2/\eps_0} \kappa (X^0_s) ds\right]  +
 \sup_{x\in \R^d} \left(\E_x \left[ \left(\int_0^{\tau_{B(x, r)}} \kappa (X^\kappa_s) ds\right)^2\right]
\P^\kappa_x \left( \tau_{B(x, r)} > kr^2/\eps_0\right) \right)^{1/2}    \nonumber \\
&\leq & \sup_{x\in \R^d} \E_x \left[ \int_0^{kr^2/\eps_0} \kappa (X^0_s) ds\right]  + \eps_0
 \sup_{x\in \R^d} \left(\E_x \left[ \left(\int_0^{\tau_{B(x, r)}} \kappa (X^\kappa_s) ds\right)^2\right]
  \right)^{1/2}   .
\label{e:3.4}
\end{eqnarray}
In view of \eqref{e:2.1} and the Aronson's estimate \eqref{e:2.7} for $\bar X^0$,
a function $f\in \K_d$ if and only if
\begin{equation} \label{e:3.5}
\lim_{t\to 0} \sup_{x\in \R^d} \bar \E^0_x \int_0^t |f(\bar X^0_s)| ds =0.
\end{equation}
Let $t_0>0 $ be as in \eqref{e:2.9}. For $0< r\leq (t_0 \eps_0/k)^{1/2}$,
\begin{eqnarray}
&& \sup_{x\in \R^d} \E_x \left[ \int_0^{kr^2/\eps_0} \kappa (X^0_s) ds\right]
= \sup_{x\in \R^d} \bar \E_x^0 \left[ M_{kr^2/\eps_0} \int_0^{kr^2/\eps_0} \kappa (\bar X^0_s) ds\right]
\nonumber \\
&\leq & \sup_{x\in \R^d} \left(  \bar \E_x^0 [ M_{t_0}^2] \, \bar \E_x^0 \left[ \left(\int_0^{kr^2/\eps_0} \kappa (\bar X^0_s) ds\right)^2 \right] \right)^{1/2} \nonumber  \\
&\leq & \left(  2C_0 \, \sup_{x\in \R^d} \bar \E_x^0 \left[ \int_0^{kr^2/\eps_0} \kappa (\bar X^0_s) \left(
\int_s^{kr^2/\eps_0} \kappa (\bar X_r^0) dr \right) ds \right] \right)^{1/2} \nonumber  \\
&\leq & \left(  2C_0 \, \sup_{x\in \R^d} \bar \E_x^0 \left[ \int_0^{kr^2/\eps_0} \kappa (\bar X^0_s)
\bar \E^0_{{\bar X^0_s}}
\left[
\int_0^{kr^2/\eps_0} \kappa (\bar X_r^0) dr \right] ds \right] \right)^{1/2} \nonumber  \\
&\leq & (2C_0)^{1/2} \sup_{x\in \R^d} \bar \E_x^0 \left[ \int_0^{kr^2/\eps_0} \kappa (\bar X^0_s) ds \right],
\label{e:3.6}
\end{eqnarray}
which goes to $0$ as $r\to 0$ by \eqref{e:3.5} as $\kappa \in \K_d$.
 Similarly, by the Markov property of $X^\kappa$,
\begin{eqnarray*}
&& \sup_{x\in \R^d} \E_x \left[ \left(\int_0^{\tau_{B(x, r)}} \kappa (X^\kappa_s) ds\right)^2\right] \\
&=& 2 \sup_{x\in \R^d} \E_x \left[ \int_0^{\tau_{B(x, r)}} \kappa (X^\kappa_s)
\left( \int_s^{\tau_{B(x, r)}} \kappa (X^\kappa_r) dr \right) ds \right] \\
&=& 2 \sup_{x\in \R^d} \E_x \left[ \int_0^{\tau_{B(x, r)}} \kappa (X^\kappa_s)
\E_{X^\kappa_s} \left[\int_0^{\tau_{B(x, r)}} \kappa (X^\kappa_r) dr \right] ds \right] \\
&\leq & 2 \left( \sup_{x\in \R^d} \E_x \left[ \int_0^{\tau_{B(x, r)}} \kappa (X^\kappa_s) ds \right] \right)^2
 \leq  2.
\end{eqnarray*}
The last inequality is due to the fact that by Proposition \ref{P:2.3},
$$
\E_x \left[ \int_0^{\tau_{B(x, r)}} \kappa (X^\kappa_s) ds \right]
\leq\E_x \left[ \int_0^{\zeta^\kappa } \kappa (X^\kappa_s) ds \right] =\P_x (\zeta^\kappa <\infty) \leq 1
\quad \hbox{for every } x\in \R^d.
$$
Hence we have by \eqref{e:3.4}-\eqref{e:3.6},
$\lim_{r\to 0} \sup_{x\in \R^d} \E_x \int_0^{\tau_{B(x, r)}} \kappa (X^\kappa_s) ds
\leq 2\eps_0 .$
 Since $\eps_0  >0$ is arbitrary, we get the the desired conclusion.
\qed

\begin{defn}\label{D:3.3} \rm
 Let $X =\{X_t, t\geq 0; \P^\kappa_x, x\in \R^d\}$ be a strong Markov process on $\R^d$ with infinitesimal generator
$\LL$.

\bee
\item[(i)] Let $A\subset \R^d$ be a measurable set.  A point $z\in \R^d$ is said to be {\it regular  for} $A$ with respect to $X $ (or equivalently, with respect to $\LL $)
if $\P_z (\sigma_A =0)=1$, where $\sigma_A:=\inf\{t>0: X_t\ \in A\}$.

\item[(ii)] An open subset $D$ of $\R^d$ is said to be {\it regular}  with respect to $X$ (or equivalently, with respect to $\LL$)
if every boundary point of $D$ is regular for $D^c$.
\eee

A celebrated result  of Littman, Stampacchia and Weinberger \cite{Litt}  states
that for an open set $D\subset \R^d$, a boundary point $z\in \partial D$ is regular for $D^c$ with respect to the Laplace  operator $\Delta$ on $\R^d$
if and only if it is regular for $D^c$ with respect to the divergence form operator $\frac12 \nabla (A(x) \nabla)$ on $\R^d$
 that is uniformly elliptic and bounded.
It follows from \cite[Proposition 3.7]{CZh1} that if  $z\in \partial D$ is regular for $D^c$ with respect to   $\Delta$, then it is regular with respect to
the operator $\LL^0$ of \eqref{e:2.3}.
 \end{defn}

 The following is the  main result of this section.

\begin{thm}\label{T:3.4}  Let $D$ be a bounded connected open subset of $\R^d$  that is  regular with respect to
the Laplace operator $\Delta$ (or equivalently, with respect to Brownian motion) on $\R^d$.
Define
$\tau_D:=\inf\{t>0: X_t\notin D\}$.
 For every  $\varphi \in C_b (D^c)$,
$ u(x):= {\mathbb E}_x \left[ \varphi(X_{\tau_D})\right]$
   is the unique weak solution of $\LL u =0$ in $D$ that is continuous on $\R^d$
with  $u(x)=\varphi (x)$ on $D^c$.
\end{thm}

\pf Let $\tau:=\inf\{t>0: X_t \not= X_{t-}\}$ be the first jumping time of $X$.
We first show that
\begin{equation}\label{e:3.7}
\P_x (\tau=\tau_D \hbox{ and } X_{\tau-} \in \partial D)=0
\quad \hbox{for every } x \in D.
\end{equation}
Denote by $\zeta^\kappa$ the lifetime of the Hunt process $X^\kappa$ associated with $\LL^0-\kappa$.
 Let $\{D_k; k\geq 1 \}$ be an increasing sequence of relatively compact subsets of $D$ that increases to $D$.
Since $\tau^\kappa_{D_k}$ is strictly increasing to $\zeta^\kappa$  on $\{\tau^\kappa _D = \zeta^\kappa \hbox{ and } X^\kappa_{\zeta^\kappa-} \in \partial D\}$, by the quasi-left continuity of $X^\kappa$, we have
for every $x\in D$,
$$
\P_x \left(\tau^\kappa _D = \zeta^\kappa \hbox{ and } X^\kappa_{\zeta^\kappa-} \in \partial D \right)
\leq \P_x \left(X^\kappa_{\zeta^\kappa-} =\lim_{k\to \infty} X^\kappa_{\tau^\kappa_{D_k}}=X^\kappa_{\zeta^\kappa} =\partial \hbox{ and }
X^\kappa_{\zeta^\kappa-} \in \partial D \right) =0.
$$
It follows then
\begin{eqnarray*}
 \P_x (\tau=\tau_D \hbox{ and } X_{\tau-} \in \partial D)
= \P_x \left(\tau^\kappa _D = \zeta^\kappa \hbox{ and } X^\kappa_{\zeta^\kappa-} \in \partial D \right)=0
\quad \hbox{for every } x\in D.
\end{eqnarray*}
This proves the claim \eqref{e:3.7}

\smallskip

(i) Let $u(x)=\E_x \left[ \varphi (X_{\tau_D})\right]$. For every $x\in D$, we have by   \eqref{e:3.1},
\begin{eqnarray*}
u(x) &=& \E_x [ \varphi (X_{\tau_D}); \tau_D <\tau] + \E_x [ \varphi (X_{\tau_D}); \tau=\tau_D
\hbox{ and } X_{\tau-} \in \partial D] \\
&& + \E_x [ \varphi (X_{\tau_D}); \tau \leq \tau_D  \hbox{ and } X_{\tau-} \in D ]   \\
&=& \E_x [ \varphi (X^\kappa_{\tau_D})]  + \E_x \left[ \int_{\R^d} \E_y \left[\varphi (X_{\tau_D})\right] \nu (X^{\kappa, D}_{\zeta^\kappa-}, dy)\right] \\
&=& \E_x [ \varphi (X^\kappa_{\tau_D})]  + \E_x \int_0^{\tau^\kappa_D} f(X^\kappa_s) ds  \\
&=:& u_1(x)+u_2(x).
\end{eqnarray*}
where $f(x):=\kappa (x) \int_{\R^d} \E_y \left[\varphi (X_{\tau_D})\right] \nu (x, dy)$, and Proposition \ref{P:2.3} is used in the third equality.
Note that $|f(x)| \leq \| u\|_\infty \kappa (x)$.
  By the proof of \cite[Theorem 3.2]{CZh2} (from line -5 on p.307 to  (3.9) on p.308 there),
   $u_2:= \E_x \int_0^{\tau^\kappa_D} f(X^\kappa_s) ds$ is in $W^{1,2}_0(D)$ and it is the unique weak solution to $(\LL^0 -\kappa )u_2 =-f$ in $D$.
Moreover, by the proof of \cite[Lemmas 5.7 and 5.8]{CZh1}, $u_2$ is continuous on $\overline D$ with $u_2=0$ on $\partial D$.
From \cite[Theorem 1.1]{CZh1}, we know
$u_1(x)= \E_x [ \varphi (X^\kappa_{\tau_D})]$ is the unique weak solution of
$(\LL^0 - \kappa) u =0$ that is continuous on $\overline D$ with $u_1=\varphi$ on $\partial D$.
We conclude that $u=u_1+u_2$ is in $W^{1,2}_{\loc }(D)\cap C(\overline D)$ with $u=\varphi$ on $\partial D$
and $u$ is a weak solution of $(\LL^0 -\kappa) u=-f=-\int_{\R^d} u(y) J (x, dy)$ in $D$.
The latter is equivalent to
$$
\LL u(x)=\LL^0 u (x) + \int_{\R^d} ( u(y)-u(x)) J(x, dy) =0 \quad \hbox{in } D.
$$
Clearly, $u(x)=\varphi (x)$ for $x\in {\overline D}^c$.
This establishes the existence part of the theorem.

\smallskip

(ii) We next show the uniqueness.
Suppose that $u$ is the unique weak solution of $\LL u =0$ in $D$ that is continuous on $\R^d$  with  $u(x)=\varphi (x)$ on $D^c$. We will show that $u(x)=\E_x [ \varphi (X_{\tau_D})]$ for $x\in D$.
 By Lemma \ref{L:3.2}, there is $r_1>0$ so that
\begin{equation}\label{e:3.8}
 \sup_{x\in \R^d} \E_x^\kappa \int_0^{\tau_{B(x, r_1)}} \kappa (X^\kappa_s) ds \leq 1/4.
\end{equation}
 For every $x_0\in D$, let $B(x_0)= B(x_0, \delta_D(x) \wedge r_1)$, where $\delta_D(x):=\inf\{|y-x|: \, y\in D^c\}$. Let $v(x):= u(x)-\E_x \left[ u(X_{\tau_{B(x_0)}}) \right]$. Then we know from (i) that $v$ is a weak solution of $\LL v=0$ in $B(x_0)$, $v\in C(\R^d)$
with $v=0$ on $B(x_0)^c$.  For    $f(x):=\int_{\R^d} v(y) J(x, dy)=\kappa (x) \int_{\R^d} v(y) \nu (x, dy)$,
 $v$ is a weak solution to
$$
 ( \kappa - \LL^0)v=   f \quad \hbox{ in } B(x_0)  \qquad \hbox{with} \quad   v=0  \hbox{ on }  B(x_0)^c.
$$
   Observe that
 $|f(x)| \leq \|v\|_\infty \kappa (x)$ and so  $ f \in \K_d$.
 By the proof of \cite[Theorem 3.2]{CZh2}  (from line -5 on p.307 to  (3.9) on p.308 there) again, we have
$  v(x)= \E_x \int_0^{\tau_{B(x_0)}} f (X^\kappa_s) ds$ for $ x\in B(x_0).$
 Hence   by Lemma \ref{L:3.1} and  \eqref{e:3.8},
\begin{eqnarray*}
 \|v\|_\infty  \leq    \|v\|_\infty  \,  \max_{x\in B(x_0)} \E_x
\int_0^{\tau_{B(x_0)}} \kappa (X^\kappa_s) ds
\leq  \frac14  \| v\|_\infty .
\end{eqnarray*}
Thus $\|v \|_\infty =0$. This proves that
\begin{equation}\label{e:3.9}
u(x) = \E_x \left[ u(X_{\tau_{B(x_0)}})\right] \quad \hbox{for every } x\in B(x_0).
\end{equation}
Take an increasing sequence $\{D_k; k\geq 1\}$ of relative compact subdoamins of $D$ that increases to $D$.
Since \eqref{e:3.9} holds for every $x_0\in D$, using the strong Markov property of $X$
we can deduce by the same argument as that for \cite[Theorem 2.2]{CS} that for every $k\geq 1$,
$u(x)= \E_x \left[ u(X_{\tau_{D_k}} )\right]$ for $ x\in D_k.$
Since $u\in C_b(\R^d)$, we get after letting $k\to \infty$
$$
u(x) = \E_x \left[ u(X_{\tau_D} )\right]= \E_x \left[ \varphi (X_{\tau_{D}})\right], \quad x\in D.
$$
This establishes the uniqueness.
\qed

\end{document}